\documentclass{amsart}
\input{diagxy}
\usepackage{amssymb}
\usepackage{stmaryrd}
\newtheorem{theorem}{Theorem}[section]
\newtheorem{lemma}[theorem]{Lemma}
\newtheorem{corollary}[theorem]{Corollary}
\newtheorem{proposition}[theorem]{Proposition}
\theoremstyle{definition}
\newtheorem{definition}[theorem]{Definition}
\newtheorem{example}[theorem]{Example}
\newtheorem{remark}[theorem]{Remark}
\newcommand{\abs}[1]{\lvert#1\rvert}
\newcommand{\cSet}{\mathbf{Set}}
\newcommand{\cCal}{\mathbf{Cl}}
\newcommand{\cC}{\mathbf{C}}
\newcommand{\A}{\mathsf{A}}
\newcommand{\B}{\mathsf{B}}
\newcommand{\C}{\mathsf{C}}
\renewcommand{\l}{\mathcal{L}}
\newcommand{\Chu}{\mathbf{Chu}}
\newcommand{\lp}{\left(}
\newcommand{\rp}{\right)}
\newcommand{\lac}{\left\{}
\newcommand{\rac}{\right\}}
\def\s{\Sigma}
\def\t{\times}
\def\aut{\mathsf{Aut}}
\def\uni{\mathsf{U}}
\def\Sep{\mathcal{S}}
\def\sep{\mathsf{S}}
\def\h{\mathcal H}
\def\pro{\mathsf{P}}
\def\p{\perp}
\def\pr{\pi^{-1}}
\def\ot{\otimes}
\newcommand{\otchu}{\mathop{\underset{\mathrm{Chu}}{\otimes}}}
\newcommand{\otaerts}{{\underset{\mbox{\rm\tiny Aerts}}{\otimes}}}
\begin{document}
\title{Weak tensor products of complete atomistic lattices}
\author{Boris Ischi}
\address{Boris Ischi, Laboratoire de
Physique des Solides, Universit\'e Paris-Sud, B\^atiment 510,
91405 Orsay, France}\curraddr{Coll\`ege de Candolle, 5 rue
d'Italie, 1204 Geneva, Switzerland}
\email{boris.ischi@edu.ge.ch}
\thanks{Supported by the Swiss National Science Foundation.}
\begin{abstract}
Given two complete atomistic lattices $\l_1$ and $\l_2$, we define
a set $\sep(\l_1,\l_2)$ of complete atomistic lattices by means of
three axioms (natural regarding the description of {\it separated}
quantum compound systems), or in terms of a universal property
with respect to a given class of bimorphisms. We call the elements
of $\sep(\l_1,\l_2)$ weak tensor products of $\l_1$ and $\l_2$. We
prove that $\sep(\l_1,\l_2)$ is a complete lattice. We compare the
bottom element $\l_1\varowedge\l_2$ with the {\it separated product}
of Aerts and with the {\it box product} of Gr\"atzer and Wehrung.
Similarly, we compare the top element $\l_1\varovee\l_2$ with the
tensor products of Fraser, Chu and Shmuely. With some additional
hypotheses on $\l_1$ and $\l_2$ (true for instance if $\l_1$ and
$\l_2$ are moreover irreducible, orthocomplemented and with the
covering property), we characterize the automorphisms of weak
tensor products in terms of those of $\l_1$ and $\l_2$.
\end{abstract}
\maketitle
\section{Introduction}

In quantum logic, on associates to a physical system $S$ a couple
$(\s_S,\l_S\subseteq 2^{\s_S})$, where $\s_S$ represents the set
of all possible states of $S$ and $\l_S$, the set of experimental
propositions concerning $S$: A proposition represented by some
$a\in\l_S$ is true with probability $1$ if and only if the state
of $S$ lies in $a$ \cite{Birkhoff/Neumann:1936}. It is usually
assumed that $\l_S$ is a {\it simple closure space} on $\s_S$ - in
other words, that $\l_S$ is closed under arbitrary intersections,
and contains $\emptyset$, $\s_S$, and all singletons of $\s_S$. In
particular, then, $\l_S$ is a complete atomistic lattice. Note
that any complete atomistic lattice is isomorphic to a simple
closure space on its set of atoms.

If $S$ is a compound system consisting of two {\it separated}
quantum systems $S_1$ and $S_2$ (for instance two electrons
prepared in two different rooms of the lab), then we have
\cite{Aerts:1982,IschiDenver:2004}
\renewcommand\theequation{P1}
\begin{equation}
\s_S=\s_{S_1}\t\s_{S_2} .
\end{equation}
In particular, simultaneous experiments on both systems can be
performed and any experiment done on one system does not alter the
state of the other system. Therefore, if $P_1$ is a proposition
concerning $S_1$ represented by some $a_1\in\l_{S_1}$ and $P_2$ is
a proposition concerning $S_2$ represented by some
$a_2\in\l_{S_2}$, then the proposition $P_1$ OR $P_2$ concerning
the compound system is true with probability $1$ if and only if
the state $p_1$ of $S_1$ lies in $a_1$ or the state $p_2$ of $S_2$
lies in $a_2$. In other words \cite{Aerts:1982,IschiDenver:2004}
\renewcommand\theequation{P2}
\begin{equation}
a_1\t\s_{S_2}\cup\s_{S_1}\t a_2\in\l_{S}\, ,\
\mbox{for all}\  a_1\in\l_{S_1}\ \mbox{and}\ a_2\in\l_{S_2} .
\end{equation}
Moreover, in addition to Axioms P1 and P2 above, we postulate that
\renewcommand\theequation{P3}
\begin{equation}
\begin{split}\mbox{For all}\
p_i\in\s_{S_i}\ \mathrm{and}\ \ A_i\subseteq\s_{S_i}\, ,\
&[p_1\t A_2\in\l_S\Rightarrow A_2\in\l_2]\ \mbox{and}\\
&[A_1\t p_2\in\l_S\Rightarrow A_1\in\l_1]\,. \end{split}
\end{equation}
See \cite{IschiDenver:2004} for detailed physical justifications
of Axioms P1-P3. We define $\sep(\l_{S_1},\l_{S_2})$ to be the
set of all simple closure spaces on $\s_{S_1}\t\s_{S_2}$ for which
Axioms P2 and P3 hold.

The rest of the paper is organized as follows. In Section
\ref{SectionMainDefinitions}, we fix some basic terminology and
notation. We define the set $\sep\equiv\sep(\l_1,\cdots,\l_n)$ of
$n$-fold weak tensor products. We prove that $\sep$ is a complete
lattice, the bottom and top elements of which are denoted by
$\varowedge_i\l_i$ and $\varovee_i\l_i$ respectively. In Section
\ref{SectionComparison} we compare $\varowedge$ and $\varovee$ with
several other tensor products of lattices. An equivalent
definition of $\sep$, in terms of a universal property with
respect to a given class of multimorphisms (arbitrary joins are
preserved) is given in Section \ref{SectionEquivalentDefinition}.
In Section \ref{SectionCentralElements} we prove some basic
relations between the central elements of $\l_i$'s and those of
$\varovee_\alpha\l_\alpha$ and $\varowedge_\alpha\l_\alpha$.

Let $\l\in\sep$ and let $u:\l\rightarrow\l$ preserve arbitrary
joins and send atoms to atoms. In Section
\ref{SectionAutomorphisms} we prove, under some hypotheses on the
image of $u$ and on each $\l_i$ (true for instance if $\l_i$ is
moreover irreducible, orthocomplemented and with the covering
property), that there exists a permutation $f$ of
$\{1,2,\cdots,n\}$, and join-preserving maps
$v_i:\l_i\rightarrow\l_{f(i)}$ sending atoms to atoms, such that
on atoms $u=f\circ (v_1\t\cdots\t v_n)$. A time evolution can be
modelled by a map preserving arbitrary joins and sending atoms to
atoms \cite{Faure/Moore/Piron:1995}. Hence, in the physical
interpretation, this result shows that separated quantum systems
remain separated only if they do not interact.
\section{Main definitions}\label{SectionMainDefinitions}

In this section we give our main definitions. We start with some
background material and basic notations used in the sequel.

\begin{definition}
Let $\s$ be a non-empty set and $\omega$ a set of subsets of $\s$. 
We write $\bigcap\omega$ (respectively $\bigcup\omega$) instead of $\bigcap_{a\in\omega}a$
(respectively $\bigcup_{a\in\omega}a$).
By a {\it simple closure space} $\l$
on $\s$ we mean a set of subsets of $\s$, ordered by
set-inclusion, closed under arbitrary set-intersections ({\it
i.e.}, for all $\omega\subseteq \l$, $\bigcap\omega\in\l$), and
containing $\s$, $\emptyset$, and all singletons. We denote the
bottom ($\emptyset$) and top ($\s$) elements by $0$ and $1$
respectively. For $p\in\s$, we identify $p$ with $\{p\}\in\l$.
Hence $p\cup q$ stands for $\{p,q\}$.

We denote the category of simple closure spaces with maps
preserving arbitrary joins (hence $0$) by $\cCal$, and the
sub-category of simple closure spaces on a particular (nonempty)
set $\s$, by $\cCal(\s)$. Let $2$ denote the simple closure space
isomorphic to the two-element lattice.

Let $\l,\,\l_1\in\cCal(\s)$ and $u:\l\rightarrow\l_1$ a map
sending atoms to atoms. We write $\aut(\l)$ for the group of
automorphisms of $\l$ and we also call $u$ the mapping from $\s$
to $\s_1$ induced by $u$.
\end{definition}

\begin{remark} Let $\l$ be a simple closure space on a (nonempty) set $\s$.
Then $\l$ is a complete atomistic lattice, the atoms of which
correspond to the points ({\it i.e.}, singletons) of $\s$. Note
that if $A\subseteq\s$, then
$\bigvee_\l(A)=\bigcap\{b\in\l\,;\,A\subseteq b\}$.

Conversely, let $\l$ be a complete atomistic lattice. Let $\s$
denote the set of atoms of $\l$, and, for each $a\in\l$, let
$\s[a]$ denote the set of atoms under $a$. Then
$\{\s[a]\,;\,a\in\l\}$ is a simple closure space on $\s$,
isomorphic to $\l$.

Finally, let $\l$ be a simple closure space and let
$u\in\aut(\l)$. Then, note that for all $a\in\l$,
$u(a)=\{u(p)\,;\,p\in a\}$.
\end{remark}

\begin{definition}\label{DefinitionAlex}
Let $\{\s_\alpha\}_{\alpha\in\Omega}$ be a family of nonempty
sets, $\mathbf{\s}=\prod_\alpha\s_\alpha$, $p\in\mathbf{\s}$,
$R\subseteq\mathbf{\s}$, $A\in\prod_\alpha 2^{\s_\alpha}$,
$\beta\in\Omega$, and $B\subseteq \s_\beta$. We shall make use of
the following notations:
\begin{enumerate}
\item We denote by $\pi_\beta:\mathbf{\s}\rightarrow\s_\beta$ the
$\beta-$th coordinate map, {\it i.e.}, $\pi_\beta(p)=p_\beta$.
\item We denote by $p[-,\beta]:\s_\beta\rightarrow\mathbf{\s}$ the
map that sends $q\in\s_\beta$ to the element of $\mathbf{\s}$
obtained by replacing $p$'s $\beta-$th entry by $q$.
\item We define $R_\beta[p]=\pi_\beta(p[\s_\beta,\beta]\cap
R)$. Note that $R_\beta[p]=\{q\in\s_\beta\,;\, p[q,\beta]\in R\}$.
\item We define $A[B,\beta]\in\prod_\alpha
2^{\s_\alpha}$ as $A[B,\beta]_\beta=B$ and
$A[B,\beta]_\alpha=A_\alpha$ for $\alpha\ne\beta$.
\item We write $\overline{A}:=\prod_\alpha A_\alpha$ and
$\overline{A}[B,\beta]:=\overline{A[B,\beta]}$.
\end{enumerate}
We omit the $\beta$ in $p[-,\beta]$ when no confusion can occur.
For instance, we write $p[\s_\beta]$ instead of
$p[\s_\beta,\beta]$.
\end{definition}

\begin{remark}\label{Remarkp[R[p]]=p[S]capR}
$p[R_\beta[p]]=p[\s_\beta]\cap R$.
\end{remark}

\begin{definition}\label{DefinitionPTensorProduct}
Let $\{\l_\alpha\}_{\alpha\in\Omega}$ be a family of simple
closure spaces on $\s_\alpha$. We denote by
$\sep(\l_\alpha,\alpha\in\Omega)$ the set all
$\l\in\cCal(\mathbf{\s})$ such that
\begin{enumerate}
\item[(P1)] $\mathbf{\s}=\prod_\alpha\s_\alpha$,
\item[(P2)] $\bigcup_\alpha \pr_\alpha(a_\alpha)\in\l$, for all
$a\in\prod_\alpha\l_\alpha$,
\item[(P3)] for all $p\in\mathbf{\s}$, $\beta\in\Omega$, and
$B\subseteq\s_\beta$, [$p[B,\beta]\in\l\Rightarrow B\in\l_\beta$].
\end{enumerate}
Let $T=\prod_\alpha T_\alpha$ with $T_\alpha\subseteq
\aut(\l_\alpha)$. We denote by $\Sep_T(\l_\alpha,\alpha\in\Omega)$
the set of all $\l\in\sep(\l_\alpha,\alpha\in\Omega)$ such that
\begin{enumerate}
\item[(P4)] for all $ v\in T$, there is $u\in \aut(\l)$ such that
$u(p)_\alpha=v_\alpha(p_\alpha)$ for all $p\in\mathbf{\s}$ and all
$\alpha\in\Omega$.
\end{enumerate}
We call elements of $\sep(\l_\alpha,\alpha\in\Omega)$ {\it weak
tensor products} of the family $\{\l_\alpha\}_{\alpha\in\Omega}$.
\end{definition}

\begin{remark} The $u$ in Axiom P4 is necessarily unique.
Note also that for all $T=\prod_\alpha T_\alpha$ with
$T_\alpha\subseteq\aut(\l_\alpha)$,
$\Sep_T(\l_\alpha,\alpha\in\Omega)\subseteq\sep(\l_\alpha,\alpha\in\Omega)$. The name ``weak tensor
product'' is justified by the fact that
$\sep(\l_\alpha,\alpha\in\Omega)$ and
$\Sep_T(\l_\alpha,\alpha\in\Omega)$ can be defined in terms of a
universal property with respect to a given class of multimorphisms
of $\cCal$ (see Section \ref{SectionEquivalentDefinition}).
\end{remark}

\begin{lemma}\label{LemmaAlex}
Let $\{\l_\alpha\}_{\alpha\in\Omega}$ be a family of simple
closure spaces on $\s_\alpha$, $\beta\in\Omega$, and
$\l\in\sep(\l_\alpha,\alpha\in\Omega)$.
\begin{enumerate}
\item For any $a\in\prod_\alpha\l_\alpha$, $\overline{a}\in\l$.
\item For any $b\in\l_\beta$ and every $p\in
\prod_\alpha\s_\alpha$, $p[b,\beta]\in\l$.
\item For any $B\subseteq\l_\beta$ and every $a\in
\prod_\alpha\l_\alpha$, $\overline{a}[\bigvee
B,\beta]=\bigvee_{b\in B} \overline{a}[b,\beta]$.
\end{enumerate}
\end{lemma}

\begin{proof} (1) Let $\beta\in\Omega$. Define
$\widehat{a}^\beta\in\prod_\alpha\l_\alpha$ by setting
$\widehat{a}^\beta_\alpha=a_\beta$ if $\alpha=\beta$, and
$\emptyset$ otherwise. Note that
$\pr_\beta(a_\beta)=\bigcup_\alpha
\pr_\alpha(\widehat{a}^\beta_\alpha)$. Now, by Axiom P2,
$\bigcup_\alpha \pr_\alpha(\widehat{a}^\beta_\alpha)\in\l$,
therefore $\pr_\beta(a_\beta)\in\l$. Finally
$\overline{a}=\bigcap_\alpha\pr_\alpha(a_\alpha)$, hence
$\overline{a}\in\l$.

(2) Define $a\in\prod_\alpha\l_\alpha$ as $a_\alpha=p_\alpha$ if
$\alpha\ne\beta$ and $a_\beta=b$. Then $p[b,\beta]=\prod_\alpha
a_\alpha$, hence from the first part $p[b,\beta]\in\l$.

(3) By the first part, $\overline{a}[\bigvee B,\beta]\in\l$.
Moreover, $\overline{a}[b,\beta]\subseteq \overline{a}[\bigvee
B,\beta]$ for all $b\in B$, hence $\bigvee_{b\in
B}\overline{a}[b,\beta]\subseteq \overline{a}[\bigvee B,\beta]$.
As a consequence, there is $X\subseteq\s_\beta$ such that
$\bigvee_{b\in B}\overline{a}[b,\beta]=\overline{a}[X,\beta]$ and
$X\subseteq\bigvee B$. Moreover, $\overline{a}[X,\beta]\in\l$,
therefore $p[X,\beta]\in\l$ for any $p\in \overline{a}$, whence,
by Axiom P3, it follows that $X\in\l_\beta$. Finally, $B\subseteq
X$, therefore $\bigvee B\subseteq X$. As a consequence, $X=\bigvee
B$.
\end{proof}

\begin{lemma}\label{LemmaP2P3}
Let $\{\l_\alpha\}_{\alpha\in\Omega}$ be a family of simple
closure spaces on $\s_\alpha$,
$\mathbf{\s}=\prod_\alpha\s_\alpha$, and $\l_0,\ \l$ and $\l_1$ be
simple closures spaces on $\mathbf{\s}$. Suppose that
$\l_0\subseteq\l \subseteq\l_1$.
\begin{enumerate}
\item If Axiom P2 holds in $\l_0$, then it holds also in $\l$.
\item If Axiom P3 holds in $\l_1$, then it holds also in $\l$.
\end{enumerate}
\end{lemma}

\begin{proof}
Direct from Definition \ref{DefinitionPTensorProduct}.
\end{proof}

\begin{definition}\label{DefinitionSeparatedProduct}
Let $\{\s_\alpha\}_{\alpha\in\Omega}$ be a family of nonempty sets
and $\{\l_\alpha\subseteq 2^{\s_\alpha}\}_{\alpha\in\Omega}$. Let
$\mathbf{\s}=\prod_\alpha\s_\alpha$. Then
\[\begin{split}\varowedge\{\l_\alpha\,;\,\alpha\in\Omega\}&:=\lac\bigcap \omega\,;\,
\omega\subseteq \lac\bigcup_\alpha \pr_\alpha(a_\alpha)\,
;\, a\in\mbox{$\prod$}_\alpha\l_\alpha\rac\rac\, ,\\[2mm]
\varovee\{\l_\alpha\,;\,\alpha\in\Omega\}&:=\{R\subseteq\mathbf{\s}\,;\,
R_\beta[p]\in\l_\beta,\, \forall
p\in\mathbf{\s},\,\beta\in\Omega\}\, ,
\end{split}\]
ordered by set-inclusion.
\end{definition}

\begin{lemma}\label{LemmaAertsP2FraserP3}
Let $\{\l_\alpha\}_{\alpha\in\Omega}$ be a family of simple
closure spaces on $\s_\alpha$ with
$\mathbf{\s}=\prod_\alpha\s_\alpha$, and let $\l_0$ and $\l_1$ be
simple closures spaces on $\mathbf{\s}$.
\begin{enumerate}
\item If Axiom P2 holds in $\l_0$, then
$\varowedge_\alpha\l_\alpha\subseteq\l_0$.
\item If Axiom P3 holds in $\l_1$ and if
$p[\s_\beta,\beta]\in\l_1$, for all $p\in\mathbf{\s}$ and all
$\beta\in\Omega$, then $\l_1\subseteq\varovee_\alpha\l_\alpha$.
\end{enumerate}
\end{lemma}

\begin{proof}
(1) Direct from Definitions \ref{DefinitionPTensorProduct} and
\ref{DefinitionSeparatedProduct}.

(2) Let $R\in\l_1$, $p\in\mathbf{\s}$, and $\beta\in\Omega$. By
hypothesis, $p[\s_\beta,\beta]\in\l_1$, hence
$p[\s_\beta,\beta]\cap R\in\l_1$. Now,
$p[\s_\beta,\beta]\cap R=p[R_\beta[p]]$ (see Remark
\ref{Remarkp[R[p]]=p[S]capR}). As a consequence,
$p[R_\beta[p]]\in\l_1$ therefore, by Axiom P3,
$R_\beta[p]\in\l_\beta$.
\end{proof}

\begin{proposition}\label{PropositionAssociativity}
Let $\{\s_\alpha\}_{\alpha\in\Omega}$ be a family of nonempty
sets, and let $\{\l_\alpha\subseteq
2^{\s_\alpha}\}_{\alpha\in\Omega}$. Let
$\{\Omega_\gamma\subseteq\Omega\,;\,\gamma\in\Gamma\}$ such that
$\Omega=\coprod\{\Omega_\gamma\,;\,\gamma\in\Gamma\}$. Then
\[\begin{split}
\varowedge\{\l_\alpha\,;\,\alpha\in\Omega\}&=
\varowedge_{\gamma\in\Gamma}(\varowedge_{\alpha\in\Omega_\gamma}
\l_\alpha)\, ,\\
\varovee\{\l_\alpha\,;\,\alpha\in\Omega\}&=
\varovee_{\gamma\in\Gamma}(\varovee_{\alpha\in\Omega_\gamma} \l_\alpha)\,
.
\end{split}
\]
\end{proposition}

\begin{proof} Direct from Definition
\ref{DefinitionSeparatedProduct}. \end{proof}

\begin{theorem}\label{TheoremTopandBottomElements}
Let $\{\l_\alpha\}_{\alpha\in\Omega}$ be a family of simple
closure spaces on $\s_\alpha$,
$\mathbf{\s}=\prod_\alpha\s_\alpha$, and $T=\prod_\alpha T_\alpha$
with $T_\alpha\subseteq\aut(\l_\alpha)$.
\begin{enumerate}
\item $\varowedge_\alpha\l_\alpha$ and $\varovee_\alpha\l_\alpha$ are
simple closure spaces on $\mathbf{\s}$.
\item $\varowedge_\alpha\l_\alpha$,
$\varovee_\alpha\l_\alpha\in\Sep_T(\l_\alpha,\alpha\in\Omega)$
\item
$\sep(\l_\alpha,\alpha\in\Omega)=\{\l\in\cCal(\mathbf{\s})\,;\,
\varowedge_\alpha\l_\alpha\subseteq\l\subseteq\varovee_\alpha\l_\alpha\}$
\item $\Sep_T(\l_\alpha,\alpha\in\Omega)$ and
$\sep(\l_\alpha,\alpha\in\Omega)$, ordered by set-inclusion, are
complete lattices.
\end{enumerate}
\end{theorem}

\begin{proof} (1) Obviously, $\varowedge_\alpha\l_\alpha$ and $\varovee_\alpha\l_\alpha$
contain $\emptyset$ and $\mathbf{\s}$, and by definition
$\varowedge_\alpha\l_\alpha$ is $\cap-$closed and
$\varovee_\alpha\l_\alpha$ contains all singletons of $\mathbf{\s}$.
If $\omega\subseteq\varovee_\alpha\l_\alpha$, then
$(\bigcap\omega)_\beta[p]=\bigcap\{R_\beta[p]\,;\, R\in\omega\}$.
Moreover, for all $p\in\mathbf{\s}$, $\bigcap_\alpha
\pr_\alpha(p_\alpha)=p$. As a consequence,
$\varowedge_\alpha\l_\alpha$ and $\varovee_\alpha\l_\alpha$ are simple
closure spaces on $\mathbf{\s}$.

(2) By definition, Axiom P2 holds in $\varowedge_\alpha\l_\alpha$ and
Axiom P3 holds in $\varovee_\alpha\l_\alpha$.

Let $a\in\prod_\alpha\l_\alpha$, $\beta\in\Omega$,
$p\in\mathbf{\s}$, and $R=\bigcup_\alpha\pr_\alpha(a_\alpha)$.
Then $R_\beta[p]=\s_\beta$ if $p_\alpha\in a_\alpha$ for some
$\alpha\ne\beta$, and $R_\beta[p]=a_\beta$ otherwise. As a
consequence, $\bigcup_\alpha\pr_\alpha(a_\alpha)\in
\varovee_\alpha\l_\alpha$, hence Axiom P2 holds in
$\varovee_\alpha\l_\alpha$, therefore $\varovee_\alpha\l_\alpha\in
\sep(\l_\alpha,\alpha\in\Omega)$.

By Lemma \ref{LemmaAertsP2FraserP3} part (1),
$\varowedge_\alpha\l_\alpha\subseteq \varovee_\alpha\l_\alpha$, hence by
Lemma \ref{LemmaP2P3}, Axiom P3 holds in $\varowedge_\alpha\l_\alpha$.
As a consequence, $\varowedge_\alpha\l_\alpha\in
\sep(\l_\alpha,\alpha\in\Omega)$.

Finally, let $v\in T$, $a\in\prod_\alpha\l_\alpha$,
$p\in\mathbf{\s}$, $R\subseteq\mathbf{\s}$, and $\omega\subseteq
2^{\mathbf{\s}}$. Then,
\[ v(\bigcup_\alpha \pr_\alpha(a_\alpha))=\bigcup_\alpha
\pr_\alpha(v_\alpha(a_\alpha))\,
,\]
and
\[v(R)_\beta[v(p)]=v_\beta(R_\beta[p]) .\]
Moreover $v(\bigcap\omega)=\bigcap\{v(x)\,;\,x\in\omega\}$. Hence,
$v$ is a bijection of $\varowedge_\alpha\l_\alpha$ and of
$\varovee_\alpha\l_\alpha$, and $v$ preserves arbitrary meets, hence
also arbitrary joins.

(3) Follows directly from Lemmata \ref{LemmaP2P3} and
\ref{LemmaAertsP2FraserP3}.

(4) Let $\omega\subseteq\sep(\l_\alpha,\alpha\in\Omega)$
(respectively $\omega\subseteq\Sep_T(\l_\alpha,\alpha\in\Omega)$).
Then obviously $\bigcap\omega=\{a\in\l\,,\,\forall
\l\in\omega\}\in\sep(\l_\alpha,\alpha\in\Omega)$ (respectively
$\bigcap\omega\in\Sep_T(\l_\alpha,\alpha\in\Omega)$).
\end{proof}

\begin{remark} Let $\{\l_\alpha\}_{\alpha\in\Omega}$ be a family of
simple closure spaces on $\s_\alpha$. For all $T=\prod_\alpha
T_\alpha$ with $T_\alpha\subseteq\aut(\l_\alpha)$,
$\Sep_T(\l_\alpha,\alpha\in\Omega)$ is a complete meet-sublattice
of $\sep(\l_\alpha,\alpha\in\Omega)$. Moreover,
$\varowedge_\alpha\l_\alpha$ and $\varovee_\alpha\l_\alpha$ are the bottom
and the top elements of $\sep(\l_\alpha,\alpha\in\Omega)$
respectively.\end{remark}

\begin{proposition} Let $\l_1$ and $\l_2$ be simple closure spaces. Then
$\l_1\varowedge\l_2\cong\l_2\varowedge\l_1$,
$\l_1\varovee\l_2\cong\l_2\varovee\l_1$, and
$2\varowedge\l_1\cong\l_1$.\end{proposition}

\begin{proof} Direct from Definition
\ref{DefinitionSeparatedProduct}.\end{proof}

More generally, then, we have:

\begin{proposition} Let $\{\l_\alpha\}_{\alpha\in\Omega}$ be a
family of simple closure spaces on $\s_\alpha$,
and $T=\prod_\alpha T_\alpha$ with
$T_\alpha\subseteq\aut(\l_\alpha)$. Then there is an isomorphism
$f:\Sep_T(\l_\alpha,\alpha\in\Omega)\rightarrow
\Sep_T(2,\l_\alpha,\alpha\in\Omega)$ such that for all
$\l\in\Sep_T(\l_\alpha,\alpha\in\Omega)$, $\l\cong f(\l)$.
\end{proposition}

\begin{proof} Direct from Definition
\ref{DefinitionPTensorProduct}.\end{proof}

\begin{example}
We now consider, for contrast, two examples, well known in
many-body quantum physics, where instead of Axioms P1 and P2, we
have
\begin{enumerate}
\item[(p1)] $\exists\, f:\prod_\alpha\s_\alpha\rightarrow
\mathbf{\s}$ with $f$ injective,
\item[(p2)] $f(\prod_\alpha\s_\alpha)\cap(\bigvee
f(\bigcup_\alpha \pr_\alpha(a_\alpha)))=f\lp\bigcup_\alpha
\pr_\alpha(a_\alpha)\rp$, $\forall a\in\prod_\alpha\l_\alpha$.
\end{enumerate}

If $\h$ is a complex Hilbert space, then $\s_\h$ denotes the set
of one-dimensional subspaces of $\h$ and $\pro(\h)$ stands for the
simple closure space isomorphic to the lattice of closed subspaces
of $\h$. Moreover, we write $\uni(\h)$ for the group of
automorphisms of $\pro(\h)$ induced by unitary maps on $\h$.

Let $\h_1$ and $\h_2$ be complex Hilbert spaces,
$\l_1=\pro(\h_1)$, $\l_2=\pro(\h_2)$ and $\l=\pro(\h_1\ot\h_2)$.
Then, Axioms p1 (with $\mathbf{\s}=\s_{\h_1\ot\h_2}$), p2 and p3
(replace $p[B,\beta]$ by $f(p[B,\beta])$ in Axiom P3) hold in
$\l$. Moreover, Axiom p4 (replace $u(p)_\alpha=v_\alpha(p_\alpha)$
by $f^{-1}(u(f(p)))_\alpha=v_\alpha(p_\alpha)$ in Axiom P4) holds
for $T=\uni(\h_1)\t\uni(\h_2)$. Note that $\bigvee (A_1\t A_2)=1$,
for all $A_i\subseteq\s_{\h_i}$ with $\bigvee A_i=1$.

Let $\h$ be a complex Hilbert space and
$\mathcal{F}=\bigoplus_{n\geq 0} \h^{\ot n}$ be the Fock space
(neither symmetrized, nor antisymmetri\-zed). Let
$\l=\pro(\mathcal{F})$, and for all integer $i$, let
$\l_i=\pro(\h)$. Let $n$ be an integer and consider the family
$\{\l_i\,;\,1\leq i\leq n\}$. Then Axioms p1, p2, p3 and p4 with
$T_i=\uni(\h)$ hold in $\l$ for all $n$. For all $i$, let
$\s_i=\s_\h$. Note that in that case, for all $n$ we have
$\bigvee(\prod_{i=1}^n\s_i)\ne 1$.
\end{example}
\section{Comparison with other tensor
products}\label{SectionComparison}

In this section, we compare the bottom element $\l_1\varowedge\l_2$
with the {\it separated product} of Aerts and with the {\it box
product} $\l_1\square\l_2$ of Gr\"atzer and Wehrung. On the other
hand, we compare the top element $\l_1\varovee\l_2$ with the
semilattice tensor product of Fraser and with the tensor products
of Chu and Shmuely.
\subsection{The separated product}

\begin{definition}
A lattice $\l$ with $0$ and $1$ is {\it orthocomplemented} if
there is a unary operation $\,^\p:\l\rightarrow\l$ such that for
all $a,\,b\in\l$, $(a^\p)^\p=a$, $a\leq b$ implies $b^\p\leq
a^\p$, and $a\bigwedge a^\p=0$.

Let $\l$ be an orthocomplemented simple closure space on $\s$.
Then, for $p,\, q\in\s$, we write $p\p q$ if $p\in q^\p$, where
$q^\p$ stands for $\{q\}^\p$.
\end{definition}

\begin{remark}
Note that the binary relation $\p$ on $\s$ is symmetric,
anti-reflexive and {\it separating}, {\it i.e.}, for all $p\ne
q\in\s$, there is $r\in\s$ such that $p\p r$ and $q\not\p r$.

Conversely, let $\s$ be a set and $\p$ a symmetric, anti-reflexive
and separating binary relation on $\s$. Then
$\l=\{A\subseteq\s\,;\, A^{\p\p}=A\}$ is an orthocomplemented
simple closure space on $\s$.
\end{remark}

\begin{definition}[D. Aerts, \cite{Aerts:1982}]\label{DefinitionOrthAerts}
Let $\{\l_\alpha\}_{\alpha\in\Omega}$ be a family of
orthocomplemented simple closure spaces on $\s_\alpha$,
$\mathbf{\s}=\prod_\alpha\s_\alpha$, and let
$p,\,q\in\mathbf{\s}$. Denote by $\#$ the binary relation on
$\mathbf{\s}$ defined by $p\# q$ if and only if there is
$\beta\in\Omega$ such that $p_\beta\p_\beta q_\beta$. Then
\[\otaerts_\alpha\l_\alpha:=\{ R\subseteq \mathbf{\s};\
R^{\#\#}=R\} .\]
\end{definition}

\begin{lemma}
Let $\{\l_\alpha\}_{\alpha\in\Omega}$ be a family of
orthocomplemented simple closure spaces on $\s_\alpha$. Then
\[\otaerts_\alpha\l_\alpha=\varowedge_\alpha\l_\alpha\, ,\]
and $\otaerts_\alpha\l_\alpha$ is orthocomplemented.\end{lemma}

\begin{proof}
Let $\mathbf{\s}=\prod_\alpha\s_\alpha$. Obviously, $\#$ is
symmetric and anti-reflexive. Since $\l_\alpha$ is
orthocomplemented, $\p_\alpha$ is separating. Therefore, it
follows directly from Definition \ref{DefinitionOrthAerts} that
$\#$ is also separating. As a consequence,
$\otaerts_\alpha\l_\alpha$ is an orthocomplemented simple closure
space on $\mathbf{\s}$. Moreover, coatoms are given by
$p^\#=\bigcup_\alpha \pr_\alpha(p_\alpha^{\p_\alpha})$. As a
consequence,
$\otaerts_\alpha\l_\alpha\subseteq\varowedge_\alpha\l_\alpha$.

Let $a\in\prod_\alpha\l_\alpha$. Denote the set of coatoms of
$\l_\alpha$ above $a_\alpha$ by $\s'[a_\alpha]$. Then
$a_\alpha=\bigcap_\alpha\s'[a_\alpha]$. Moreover,
\[\bigcap\{\bigcup_\alpha
\pr_\alpha(x_\alpha)\,;\,x\in\mbox{$\prod$}_\alpha\s'[a_\alpha]\}=
\bigcup_\alpha \pr_\alpha(a_\alpha) .\]
Hence $\bigcup_\alpha\pr_\alpha(a_\alpha)\in
\otaerts_\alpha\l_\alpha$. Another way to see this is to compute
\[\lp \bigcup_\alpha\pr_\alpha(a_\alpha)\rp^\#=\bigcap_\alpha\lp
\pr_\alpha(a_\alpha)\rp^\#=\bigcap_\alpha\pr_\alpha(a_\alpha^{\p_\alpha})
=\prod_\alpha a_\alpha^{\p_\alpha} \, ,\]
whence
\[\lp \bigcup_\alpha\pr_\alpha(a_\alpha)\rp^{\#\#}=
\bigcup_\alpha\pr_\alpha(a_\alpha^{\p_\alpha\p_\alpha})=
\bigcup_\alpha\pr_\alpha(a_\alpha)
 .\]

As a consequence,
$\varowedge_\alpha\l_\alpha\subseteq\otaerts_\alpha\l_\alpha$.
\end{proof}

\begin{remark} The symbol $\varowedge$ was originally used by
Aerts.\end{remark}
\subsection{The box product}

\begin{definition}
Let $\l$ be a lattice and $a\in\l$. We denote by $a\downarrow$ the
set $\{x\in\l\,;\, x\leq a\}$ and by $L$ the set
$\{a\downarrow\,;\, a\in\l\}$ ordered by set-inclusion (note that
$L\subseteq 2^{\l}$).
\end{definition}

\begin{definition}[G. Gr\"atzer, F. Wehrung,
\cite{Graetzer/Wehrung:1999}]\label{DefinitionWehrung} Let $\l_1$
and $\l_2$ be lattices and $(a,b)\in\l_1\t\l_2$. Define
\[a\square b=(a\downarrow\t\l_2)\cup(\l_1\t b\downarrow) .\]
The {\it box product} $\l_1\square\l_2$ is defined as the set of
all finite intersections of the form $\bigcap \{a_i\square
b_i\,;\, i\leq n\}$ where $(a_i,b_i)\in\l_1\t\l_2$ for all $i\leq
n$, ordered by set-inclusion.
\end{definition}

\begin{remark}
Obviously, $\l_1\square\l_2$ is a meet-sublattice of
$2^{\l_1\t\l_2}$. In fact, $\l_1\square\l_2$ is a lattice (see
\cite{Graetzer/Wehrung:1999}, Proposition 2.9). Note that if
$\l_1$ and $\l_2$ have top elements, then
$\l_1\boxtimes\l_2=\l_1\square\l_2$ where $\l_1\boxtimes\l_2$ is
the {\it lattice tensor product} (see
\cite{Graetzer/Wehrung:1999}).
\end{remark}

\begin{definition}\label{DefinitionAertsFinite}
Let $\s_1$ and $\s_2$ be nonempty sets, $\l_1\subseteq 2^{\s_1}$
and $\l_2\subseteq 2^{\s_2}$. Define $\l_1\varowedge_n\l_2$ as in
Definition \ref{DefinitionSeparatedProduct} but taking only finite
intersections. If $\l_1$ and $\l_2$ are atomistic lattices, then
define $\l_1\varowedge_n\l_2$ as $l_1\varowedge_n l_2$, where
$l_i=\{\s[a_i]\,;\,a_i\in\l_i\}\subseteq 2^{\s_i}$ and $\s[a_i]$
denotes the set of atoms under $a_i$ ($\s_i:=\s[1_i]$).
\end{definition}

\begin{proposition}
Let $\l_1$ and $\l_2$ be lattices. Then
$\l_1\square\l_2=L_1\varowedge_n L_2$.\end{proposition}

\begin{proof} Direct from Definitions \ref{DefinitionWehrung} and
\ref{DefinitionAertsFinite}.\end{proof}

\begin{theorem} For atomistic lattices, $\l_1\varowedge_{n}\l_2\cong\l_1\square\l_2$.
\end{theorem}

\begin{proof} Define $f:\l_1\varowedge_n\l_2\rightarrow\l_1\square\l_2$ as
\[f(\s[a_1]\t \s_2\cup \s_1\t \s[a_2])=a_1\square
a_2\, ,\]
and set $f(\bigcap x_i):=\bigcap f(x_i)$, where
each $x_i$ has the form $\s[a_1]\t \s_2\cup \s_1\t \s[a_2]$.
Obviously, $f$ is bijective and preserves meets, hence also joins.
\end{proof}
\subsection{The tensor products of Chu and Shmuely}

We now compare $\varovee$ with three other tensor products of lattices
appearing in the literature.

\begin{definition}[P. H. Chu \cite{Barr:handbook}]
The category $\Chu(\cSet,2)$ has as objects triples $\A=(A,r,X)$
where $A$ and $X$ are sets and $r$ is a map from $A\t X$ to $2$
(where $2=\{0,1\}$). Arrows are pairs of maps
$(F,G):\A\rightarrow(B,s,Y)$ with $F:A\rightarrow B$ and
$G:Y\rightarrow X$ such that $s(F(a),y)=r(a,G(y))$ for all $a\in
A$ and $y\in Y$. There is a functor $\,^\p$ and a bifunctor
$\otchu$ defined on objects as $\A^\p=(X,\check{r},A)$ with
$\check{r}(x,a)=r(a,x)$, and $\A\otchu \B=(A\t
B,t,\Chu(\A,\B^\p))$ with
$t((a,b),(F,G))=r(a,G(b))=\check{s}(F(a),b)$. \end{definition}

\begin{remark} A simple closure space $\l\in\cCal(\s)$ can be
identified with the Chu space $(\s,r,\l)$, which we also denote by
$\l$, where $r(p,a)=1$ if and only if $p\in a$.\end{remark}

\begin{proposition}
Let $\l_0$ and $\l_1$ be simple closure spaces on $\s_0$ and
$\s_1$ respectively. Then there is a one to one correspondence
between invertible arrows in $\Chu(\l_0,\l_1)$ and isomorphisms
between $\l_0$ and $\l_1$.
\end{proposition}

\begin{proof}
Let $f:\l_0\rightarrow\l_1$ be an isomorphism. Then
$(f,f^{-1})\in\Chu(\l_0,\l_1)$. Let $(F,G)$ be an invertible arrow
in $\Chu(\l_0,\l_1)$. Then $G:\l_1\rightarrow\l_0$ is bijective,
and since $G=F^{-1}$, $G$ preserves arbitrary meets, hence also
arbitrary joins.\end{proof}

\begin{definition}Let $\l_1$ and $\l_2$ be posets. A Galois connection
(more precisely an adjunction) between $\l_1$ and $\l_2$ is a pair
$(f,g)$ of order-preserving maps with $f:\l_1\rightarrow\l_2$ and
$g:\l_2\rightarrow\l_1$, such that for all $a\in\l_1$ and
$b\in\l_2$, $f(a)\leq b\Leftrightarrow a\leq g(b)$. Note that if
$\l_1$ and $\l_2$ are complete lattices, then there is a one to
one correspondence between maps from $\l_1$ to $\l_2$ preserving
arbitrary joins and adjunctions between $\l_1$ and $\l_2$.
\end{definition}

The following is standard; we include the easy proof for
completeness.

\begin{lemma}\label{LemmaJoin-preserving}
Let $\l_0$ and $\l_1$ be complete join-semilattices and $(f,g)$ a
pair of maps with $f:\l_0\rightarrow\l_1$ and
$g:\l_1\rightarrow\l_0$.
\begin{enumerate}
\item If $(f,g)$ forms a Galois connection, then $f$ preserves
arbitrary joins.

\item Suppose that $\l_0$ is a simple closure space on $\s_0$. For
$b\in\l_1$, let $F^{-1}(b)$ denote the set
$\{p\in\s_0\,;\,f(p)\leq b\}$. Then $f$ preserves arbitrary joins
$\Leftrightarrow$ $f(a)=\bigvee \{f(p)\,;\,p\in a\}$ for all
$a\in\l_0$ and $F^{-1}(b)\in\l_0$ for all $b\in\l_1$.
\end{enumerate}
\end{lemma}

\begin{proof} (1) Let $\omega\subseteq\l_0$. Since $f$ preserves
order, $y:=\bigvee\{f(x)\,;\,x\in\omega\}\leq f(\bigvee\omega)$.
On the other hand, for any $x\in\omega$, $f(x)\leq y$, hence
$x\leq g(y)$. As a consequence, $\bigvee\omega\leq g(y)$,
therefore $f(\bigvee\omega)\leq y$.

(2) ($\Rightarrow$) If $p\in\bigvee F^{-1}(b)$, then $f(p)\leq
\bigvee \{f(q)\,;\, q\in F^{-1}(b)\}\leq b$, hence $p\in
F^{-1}(b)$. ($\Leftarrow$) Obviously, $(f,F^{-1})$ forms a Galois
connection between $\l_0$ and $\l_1$.
\end{proof}

\begin{remark}
Let $\l$ be a poset. Then $\l^*$ denotes the dual of $\l$ (defined
by the converse order-relation) and $\leq_*$ the order relation in
$\l^*$.
\end{remark}

\begin{theorem}\label{TheoremOtf=Chu}
Let $\l_1$ and $\l_2$ be simple closure spaces on $\s_1$ and
$\s_2$ respectively. Then,
\begin{enumerate}
\item $\l_1\varovee\l_2\cong\l_1\otchu\l_2$,
\item $(\l_1\varovee\l_2)^*$ is isomorphic to the set of
join-preserving maps $\l_1\rightarrow\l_2^*$, ordered point-wise
({\it i.e.}, $f\leq g$ if and only if $f(p_1)\leq_* g(p_1)$ for
all $p_1\in\s_1$).
\end{enumerate}
\end{theorem}

\begin{proof} (1) We denote $\s_1\t\s_2$ by $\mathbf{\s}$.
Let $(F,G)\in\Chu(\l_1,\l_2^\p)$. Define
\[R_{(F,G)}:=\bigcup_{p_1\in\s_1}p_1\t F(p_1) .\]
By definition, ${R_{(F,G)}}_2[p]\in\l_2$ for all
$p\in\mathbf{\s}$, and
\[q_1\in {R_{(F,G)}}_1[p]\Leftrightarrow(q_1,p_2)\in
R_{(F,G)}\Leftrightarrow p_2\in F(q_1)\Leftrightarrow q_1\in
G(p_2)\, ,\]
thus ${R_{(F,G)}}_1[p]\in\l_1$ for all $p\in\mathbf{\s}$, hence
$R_{(F,G)}\in\l_1\varovee\l_2$.

Let $R\in\l_1\varovee\l_2$. Define $F_R:\s_1\rightarrow\l_2$ and
$G_R:\s_2\rightarrow\l_1$ as $F_R(p_1):=R_2[(p_1,\cdot)]$ and
$G_R(p_2):=R_1[(\cdot,p_2)]$. Then
\[p_2\in F_R(p_1)\Leftrightarrow (p_1,p_2)\in R\Leftrightarrow
p_1\in G_R(p_2)\, ,\]
so that $(F_R,G_R)\in\Chu(\l_1,\l_2^\p)$. Note that $(F,G)\mapsto
R_{(F,G)}$ is isotone in both directions.

Obviously,
\[F_{R_{(F,G)}}=F,\ G_{R_{(F,F)}}=G,\ \mbox{and}\
R_{(F_R,G_R)}=R .\]
As a consequence, there is an invertible arrow in
$\Chu(\l_1\varovee\l_2,\l_1\otchu\l_2)$.

(2) Let $(F,G)\in\Chu(\l_1,\l_2^\p)$. Define
$f:\l_1\rightarrow\l_2^*$ and $g:\l_2^*\rightarrow\l_1$ as
$f(a):=\bigcap F(a)$ and $g(b)=\bigcap G(b)$. Then, obviously
$(f,g)$ forms a Galois connection between $\l_1$ and $\l_2^*$.

Let $f:\l_1\rightarrow\l_2^*$ be $\bigvee-$preserving. Define
$g:\l_2\rightarrow\l_1$ as $g(b):=\bigvee\{a\in\l_1\,; f(a)\leq_*
b\}$. Then, obviously $(f,g)\in\Chu(\l_1,\l_2^\p)$. \end{proof}

\begin{remark} Note that $\varovee$ is the $\boxtimes$-tensor
product  of Golfin \cite{Golfin:handbook}. As a corollary of
Theorem \ref{TheoremOtf=Chu} part (2), $\l_1\varovee\l_2\cong
\l_1\ot\l_2$ the tensor product of Shmuely \cite{Shmuely:1974}.
\end{remark}
\subsection{The semilattice tensor product of Fraser}

We now consider tensor products as solutions to universal mapping
problems.

\begin{definition} Let $\cC$ be a concrete
category over $\cSet$. Let $\{\A_\alpha\}_{\alpha\in\Omega}$ and
$\C$ be objects of $\cC$. A {\it multimorphism} of $\cC$ is an
arrow $f\in\cSet(\prod_\alpha\A_\alpha,\C)$ such that for all
$a\in\prod_\alpha\A_\alpha$, $f(a[-,\alpha])\in\cC(\A_\alpha,\C)$
(see Definition \ref{DefinitionAlex}). If $\abs{\Omega}=2$, then
$f$ is called a {\it bimorphism}.
\end{definition}

\begin{definition}[G. Seal, \cite{Seal:2004}]\label{DefinitionTensorGavin}
Let $\cC$ be a concrete category over $\cSet$. A {\it tensor
product} in ${\bf C}$ is a bifunctor
$-\ot-:\cC\t\cC\rightarrow\cC$ such that for all objects $\A$,
$\B$ of $\cC$, there is a a bimorphism $f:
\A\t\B\rightarrow\A\ot\B$, and for all object $\C$ of $\cC$, and
all bimorphism $g:\A\t\B\rightarrow\C$, there is a unique arrow
$h\in\cC(\A\ot\B,\C)$ which makes the diagram
$$\bfig
 \ptriangle|alm|/>`>`>/[\A\t\B`\A\ot\B`\C;f`g`!h]
 \efig$$
commute.
\end{definition}

\begin{remark}
By definition, the tensor product is unique up to isomorphisms.
Consider now the category $\cC$ of join-semilattices with maps
preserving finite joins. Let $\l_1$, $\l_2$ and $\l$ be
join-semilattices. If $\l$ is a tensor product of $\l_1$ and
$\l_2$ in $\cC$, then $f(\l_1\t\l_2)$ generates $\l$. Indeed, let
$\l_0$ be a (join-semilattice) tensor product of $\l_1$ and $\l_2$
as defined by Fraser in \cite{Fraser:1976} (write $g$ the
bimorphism from $\l_1\t\l_2$ to $\l_0$). Then by definition
$\l_0$ is generated by $g(\l_1\t\l_2)$, and obviously $\l_0$ is a
tensor product of $\l_1$ and $\l_2$ in $\cC$. Therefore there is
an isomorphism $h:\l_0\rightarrow\l$ such that $f=h\circ g$.
Hence, for every $x\in\l$ we have $x=h(h^{-1}(x))=h(\bigvee
g(a_i,b_i))=\bigvee h(g(a_i,b_i))=\bigvee f(a_i,b_i)$ where
$a_i\in\l_1$ and $b_i\in\l_2$.

As a consequence, for the category of join-semilattices with maps
preserving finite joins, Definition \ref{DefinitionTensorGavin} is
equivalent to the definition of Fraser.
\end{remark}

\begin{theorem}\label{TheoremTensorF=bi-tensor}
$\varovee$  is a tensor product in $\cCal$.
\end{theorem}

\begin{proof} Let $\l_1,\,\l_2\in\cCal$. By Theorem
\ref{TheoremTopandBottomElements}, $\l_1\varovee\l_2\in\cCal$, and by
Lemma \ref{LemmaAlex}.3, the map $f:\l_1\t\l_2\rightarrow
\l_1\varovee\l_2$ sending $a=(a_1,a_2)$ to $\overline{a}=a_1\t a_2$ is
a bimorphism.

Let $\l\in\cCal$ (or a complete join-semilattice), and let
$g:\l_1\t\l_2\rightarrow \l$ be a bimorphism. Define
$h:\l_1\varovee\l_2\rightarrow \l$ by $h(R)=\bigvee\{g(p_1,p_2)\,;\,
(p_1,p_2)\in R\}$. Furthermore, for $p_1\in\s_1$ and $p_2\in\s_2$,
define $g_{p_1}:\l_2\rightarrow\l$ and $g_{p_2}:\l_1\rightarrow\l$
as $g_{p_1}(y):=g(p_1,y)$ and $g_{p_2}(x):=g(x,p_2)$ respectively.
Let $b\in\l$. Write, as in Lemma \ref{LemmaJoin-preserving},
$H^{-1}(b)=\{p\in\s_1\t\s_2\,;\, h(p)\leq b\}$,
$G_{p_1}^{-1}(b)=\{s\in\s_2\,;\, g_{p_1}(s)\leq b\}$, and
$G_{p_2}^{-1}(b)=\{r\in\s_1\,;\, g_{p_2}(r)\leq b\}$. Then,
\[ H^{-1}(b)=\bigcup_{p_1\in\s_1}p_1\t
G_{p_1}^{-1}(b)=\bigcup_{p_2\in\s_2}G_{p_2}^{-1}(b)\t p_2 .\]
Since $g$ is a bimorphism, it follows from Lemma
\ref{LemmaJoin-preserving} and Definition
\ref{DefinitionSeparatedProduct} that $H^{-1}(b)\in\l_1\varovee\l_2$.
Therefore, by Lemma \ref{LemmaJoin-preserving}, $h\in
\cCal(\l_1\varovee\l_2,\l)$. Finally, if $h'\in\cCal(\l_1\varovee\l_2,\l)$
and $h'\circ f=g$, then $h'$ equals $h$ on atoms, therefore
$h'=h$.

Let $L_1,\,L_2\in\cCal$, $f_1\in\cCal(\l_1,L_1)$, and
$f_2\in\cCal(\l_2,L_2)$. Then $g=f\circ(f_1\t
f_2):\l_1\t\l_2\rightarrow L_1\varovee L_2$ is a bimorphism. We define
$f_1\varovee f_2\in\cCal(\l_1\varovee\l_2,L_1\varovee L_2)$ to be the arrow
$h$ constructed above.
\end{proof}

\begin{remark}
Let $\cC_{com}$ denote the category of complete join-semilattices
with maps preserving arbitrary joins. Let $\l_1$ and
$\l_2\in\cCal$. Then $\l_1,\,\l_2,\,\l_1\varovee\l_2\in\cC_{com}$.
Moreover, we have proved that $\l_1\varovee\l_2$ is the tensor product
of $\l_1$ and $\l_2$ in $\cC_{com}$.
\end{remark}
\section{Equivalent definition}\label{SectionEquivalentDefinition}

In this section we give an equivalent definition of
$\sep(\l_\alpha,\alpha\in\Omega)$ and
$\Sep_T(\l_\alpha,\alpha\in\Omega)$ in terms of a universal
property with respect to a given class of bimorphisms of $\cCal$.

\begin{lemma}\label{LemmaExitsh_a:L->2}
Let $\l$ be a simple closure space on $\s$. For all $a\in\l$,
there is a unique mapping $h_a:\l\rightarrow 2$ preserving
arbitrary joins such that for all $p\in\s$, $h_a(p)=0$ if and only
if $p\in a$.
\end{lemma}

\begin{proof} Let $b\in\l$.
Define $h_a$ by $h_a(b)=\bigvee\{ h_a(p)\,;\, p\in b\}$ (hence
$h_a(0)=0$) with $h_a(p)=0$ if $p\in a$ and $h_a(p)=1$ if
$p\not\in a$. Then, by Lemma \ref{LemmaJoin-preserving}, $h$
preserves arbitrary joins.\end{proof}

\begin{lemma} Let $\l$ and $\l_\alpha$ ($\alpha\in\Omega$) be
simple closure spaces. Then
\begin{enumerate}
\item If $f:\prod_\alpha\l_\alpha\rightarrow \l$ a multimorphism
of $\cCal$ and $v\in\prod_\alpha \aut(\l_\alpha)$, then $f\circ v$
is again a multimorphism.
\item For $g\in\prod_\alpha\cCal(\l_\alpha,2)$, the mapping
$\bigcap_\alpha g_\alpha:\prod_\alpha\l_\alpha\rightarrow 2$,
defined as
\[\bigcap_\alpha g_\alpha(a):=\bigcap\{
g_\alpha(a_\alpha)\,;\,\alpha\in\Omega\}\, ,\]
is a multimorphism.
\end{enumerate}
\end{lemma}

\begin{proof}
(1) Let $w=f\circ v$, $a\in\prod_\alpha\l_\alpha$,
$\beta\in\Omega$, and $\omega\subseteq\l_\beta$. Then
\[w(a[\bigvee\omega,\beta])=f(v(a)[\bigvee\{
v_\beta(x)\,;\,x\in\omega\},\beta])=\bigvee_{x\in\omega}w(a[x,\beta])\,
,\]
since $f$ is a multimorphism.

(2) Let $g=\bigcap_\alpha g_\alpha$, $a\in\prod_\alpha\l_\alpha$,
$\beta\in\Omega$, and $b\in\l$. Write
\[
h(-):=g(a[-,\beta]):\l_\beta\rightarrow 2 .
\] 
Then $h(b)=g_\beta(b)\bigcap_{\alpha\ne\beta}g_\alpha(a_\alpha)$. Hence $h(b)=0$ if $g_\alpha(a_\alpha)=0$ for $\alpha\ne
\beta$ and $h(b)=g_\beta(b)$ otherwise. As a consequence, $h$ preserves arbitrary joins,
hence $g$ is a multimorphism of $\cCal$.
\end{proof}

\begin{definition}
Let $\l$ and $\l_\alpha$ ($\alpha\in\Omega$) be simple closure
spaces, $f:\prod_\alpha\l_\alpha\rightarrow \l$ a multimorphism of
$\cCal$ and $T_\alpha\subseteq\aut(\l_\alpha)$. Put
$T:=\prod_\alpha T_\alpha$. Then, we define
\[\begin{split}
V_T(f)&:=\{f\circ v:\prod_\alpha\l_\alpha\rightarrow\l\,;\,v\in
T\}\,
,\\
\Gamma_\cap&:=\{\bigcap_\alpha
g_\alpha:\prod_\alpha\l_\alpha\rightarrow
2\,;\,g\in\prod_\alpha\cCal(\l_\alpha,2)\} .\end{split}
\]
\end{definition}
\begin{definition}Let $\l$ be a simple closure space and
$A\subseteq\l$. We say that $A$
{\it generates} $\l$ if each element of $\l$ is the join of
elements in $A$.\end{definition}

\begin{theorem}
Let $\{\l_\alpha\}_{\alpha\in\Omega}$ be a family of simple
closure spaces on $\s_\alpha$, $\l$ a simple closure space on
$\mathbf{\s}$, and $T=\prod_\alpha T_\alpha$ with
$T_\alpha\subseteq\aut(\l_\alpha)$.
\begin{enumerate}
\item If $\l\in\Sep_T(\l_\alpha,\alpha\in\Omega)$, then there is a
multimorphism $f:\prod_\alpha\l_\alpha\rightarrow \l$ of $\cCal$
such that $f(\prod_\alpha\l_\alpha)$ generates $\l$, and for all
$g\in\Gamma_\cap$ and $w\in V_T(f)$, there is a unique
$h\in\cCal(\l,2)$ and a unique $u\in\cCal(\l,\l)$ such that the following
diagrams commute.
$$\bfig
 \ptriangle|alm|/>`>`>/[\prod_\alpha\l_\alpha`\l`2;f`g`!h]
 \efig
 \hspace{1cm}\bfig
 \ptriangle|alm|/>`>`>/[\prod_\alpha\l_\alpha`\l`\l;f`w`!u]
 \efig$$
%
\item Conversely, if there is a multimorphism of $\cCal$
$f:\prod_\alpha\l_\alpha\rightarrow \l$ satisfying all conditions
of part (1), then there is some
$\l_0\in\Sep_T(\l_\alpha,\alpha\in\Omega)$ such that
$\l\cong\l_0$.
\end{enumerate}
\end{theorem}

\begin{proof} (1) By Lemma
\ref{LemmaAlex}.3, the map $f:\prod_\alpha\l_\alpha\rightarrow\l$
sending $a$ to $\overline{a}$, is a multimorphism of $\cCal$, and
obviously, $f(\prod_\alpha\l_\alpha)$ generates $\l$. Moreover, by
Axiom P4, for all $w\in V_T(f)$ there is a unique $u$ such that
the second diagram commutes.

Let $g=\bigcap_\alpha g_\alpha\in\Gamma_\cap$. Write $G_\alpha$
for the restriction to atoms of $g_\alpha$. Let $a_\alpha:=\bigvee
G_\alpha^{-1}(0)$. Note that by Lemma \ref{LemmaJoin-preserving},
$a_\alpha=G_\alpha^{-1}(0)$. Recall that by hypothesis, $\l\in\cCal(\mathbf{\s})$ with
$\mathbf{\s}=\prod_\alpha\s_\alpha$. Define $x:=\bigcup_\alpha
\pr_\alpha(a_\alpha)$. By Axiom P2 $x\in\l$, and by Lemma
\ref{LemmaExitsh_a:L->2} there is $h_x\in\cCal(\l,2)$ such that
for all $p\in\mathbf{\s}$, $h_x(p)=0$ if and only if $p\in x$.
Hence, $h_x\circ f=g$. Let $h\in\cCal(\l,2)$ such that $h\circ
f=g$. Then on atoms $h$ equals $h_x$, therefore $h=h_x$.

(2) Let $\mathbf{\s_0}=\prod_\alpha\s_\alpha$. We denote by $F$
the mapping from $\mathbf{\s_0}$ to $\l$ induced by the
multimorphism $f$.

(2.1) For $a,\, b\in\prod_\alpha\l_\alpha\backslash\{0_\alpha\}$, we write $a\leq b$ if
and only if $a_\beta\leq b_\beta$, for all $\beta\in\Omega$ (this
is the standard product ordering on $\prod_\alpha\l_\alpha$). {\it
Claim}: $f(a)\subseteq f(b)\Rightarrow a\leq b$. As a corollary,
$f$ is injective. [{\it Proof}: Suppose that $f(a)\subseteq f(b)$
and that $a_\beta\not\leq b_\beta$ for some $\beta\in\Omega$. By
Lemma \ref{LemmaExitsh_a:L->2}, there is
$h_{b_\beta}\in\cCal(\l_\beta,2)$ such that for all
$p\in\s_\beta$, $h_{b_\beta}(p)=0$ if and only if $p\in b_\beta$.
Let $g_\beta:=h_{b_\beta}$, $g_\alpha:=h_{0_\alpha}$ for all
$\alpha\ne\beta$, and $g:=\bigcap_\alpha g_\alpha$. By definition,
$g\in \Gamma_\cap$, hence there is $h\in\cCal(\l,2)$ such that
$h\circ f=g$. As a consequence,
\[1=g_{\beta}(a_\beta)=g(a)=h(f(a))\subseteq
h(f(b))=g(b)=g_\beta(b_\beta)=0\, ,\]
a contradiction. This proves the claim.]

As a consequence, since $f(\prod_\alpha\l_\alpha)$ generates $\l$,
for all $p\in\mathbf{\s_0}$, $F(p)$ is an atom of $\l$, and the
mapping from $\mathbf{\s_0}$ to $\mathbf{\s}$ induced by $F$
(which we also denote by $F$) is bijective. Moreover, for all
$a\in\prod_\alpha\l_\alpha$,
\[f(a)=\bigvee\{F(p)\,;\,p\in \overline{a}\} .\]
Therefore, for all $a,\, b\in\prod_\alpha\l_\alpha$, we have:
$a\leq b\Rightarrow f(a)\subseteq f(b)$ (note that if $\Omega$ is
finite, then this implication follows directly from the fact that
$f$ is a multimorphism).

(2.2) Let $\l_0\subseteq 2^{\mathbf{\s_0}}$ defined as
$\l_0:=\{F^{-1}(c)\,;\,c\in\l\}$. Then by what precedes,
$\l_0\in\cCal(\mathbf{\s_0})$ and the map
$F^{-1}:\l\rightarrow\l_0$ is bijective and preserves arbitrary
meets, hence also arbitrary joins. It remains to prove that
$\l_0\in\Sep_T(\l_\alpha,\alpha\in\Omega)$, hence to check that
Axioms P2, P3 and P4 hold in $\l_0$. Below, if $g\in\Gamma_\cap$,
then $G$ denotes the map from $\mathbf{\s_0}$ to $2$ induced by
$g$.

(P2) Let $a\in\prod_\alpha\l_\alpha$, $x:=F(\bigcup_\alpha
\pr_\alpha(a_\alpha))$, and $p\in\mathbf{\s_0}$ such that
$p_\alpha\not\in a_\alpha$ for all $\alpha\in\Omega$. Suppose that
$F(p)\in \bigvee x$. From Lemma \ref{LemmaExitsh_a:L->2}, for all
$\beta\in\Omega$, there is $h_{a_\beta}\in\cCal(\l_\beta,2)$ such
that for all $p\in\s_\beta$, $h_{a_\beta}(p)=0$ if and only if
$p\in a_\beta$. Let $g=\bigcap_\alpha h_{a_\alpha}$. By
definition, $g\in \Gamma_\cap$, hence there is $h\in\cCal(\l,2)$
such that $h\circ f=g$; whence
\[\begin{split}1=G(p)=h(F(p))\subseteq h(\vee x)&=
\bigvee \{h(F(p))\,;\,p\in\bigcup_\alpha\pr_\alpha(a_\alpha)\} \\
&=\bigvee \{G(p)\,;\,p\in\bigcup_\alpha\pr_\alpha(a_\alpha)\}=0\,
,
\end{split}\]
a contradiction. As a consequence, $F^{-1}(\bigvee
F(\bigcup_\alpha \pr_\alpha(a_\alpha)))=\bigcup_\alpha
\pr_\alpha(a_\alpha)$, hence $\bigcup_\alpha
\pr_\alpha(a_\alpha)\in\l_0$.

(P3) Let $\beta\in\Omega$, $A\subseteq\s_\beta$, and
$p\in\mathbf{\s_0}$, such that $p[A,\beta]\in\l_0$, {\it i.e.},
there is $c\in\l$ such that $p[A,\beta]=F^{-1}(c)$. Let $q\in
\bigvee A$ and $\hat{p}\in\prod_\alpha\l_\alpha$ such that
$\overline{\hat{p}}=\{p\}$. Then $p[q,\beta]\in p[\bigvee
A,\beta]$, hence, since $f$ is multimorphism, we find that
\[ F(p[q,\beta])\in F(p[\bigvee A,\beta])\subseteq\bigvee F(p[\bigvee
A,\beta])=f(\hat{p}[\bigvee A,\beta])=\bigvee_{q\in A}
F(p[q,\beta])\subseteq c .\]
As a consequence, $q\in A$, therefore $A\in\l_\beta$.

(P4) Let $v\in T$ and $u\in\cCal(\l,\l)$ such that $u\circ
f=f\circ v$. Define $u_0:=F^{-1}\circ u\circ F$. Then
$u_0\in\aut(\l_0)$ and $u_0(p)_\alpha=v_\alpha(p_\alpha)$ for all
$p\in\mathbf{\s_0}$ and all $\alpha\in\Omega$.
\end{proof}
\section{Central elements}\label{SectionCentralElements}

Let $\{\l_\alpha\}_{\alpha\in\Omega}$ be a family of simple
closure spaces and $\beta\in\Omega$. In this section we prove that
if $z$ is a central element of $\l_\beta$, then $\pr_\beta(z)$ is
a central element of $\varovee_\alpha\l_\alpha$ and of
$\varowedge_\alpha\l_\alpha$. As a corollary, if
$\varovee_\alpha\l_\alpha$ or if $\varowedge_\alpha\l_\alpha$  is
irreducible, then all $\l_\alpha$'s are irreducible. We give some
sufficient conditions under which the converse result holds in
$\varowedge_\alpha\l_\alpha$.

\begin{definition}
Let $a$ and $b$ be elements of a lattice. Then $(a,b)$ is said to
be a {\it modular pair} (in symbols $(a,b)M$) if $(c\vee
a)\wedge b=c\vee(a\wedge b)$ for all $c\leq b$.
\end{definition}

\begin{lemma}\label{LemmaCentralElements}
Let $\l$ be a simple closure space on $\s$ and $z\in\l$. Then $z$
is a central element of $\l$ if and only if $z^c:=\s\backslash
z\in\l$ and $(z,z^c)M$ and $(z^c,z)M$.
\end{lemma}

\begin{proof} Direct from Theorem 4.13 ($\varepsilon$) in
\cite{Maeda/Maeda:handbook}.\end{proof}

\begin{corollary}\label{CorollaryPatula}
Let $\l$ be  a simple closure space on $\s$, $z$ a central
element, $a\subseteq z$, and $b\subseteq z^c$. Then $a\vee
b=a\cup b$.
\end{corollary}

\begin{proof} From $(z^c,z)M$ follows that
$(a\vee b)\cap z\subseteq(a\vee z^c)\cap z= a$, and
from $(z,z^c)M$ follows that $(a\vee b)\cap
z^c\subseteq(b\vee z)\cap z^c= b$.\end{proof}

\begin{theorem} Let $\{\l_\alpha\}_{\alpha\in\Omega}$ be a
family of simple closure spaces on $\s_\alpha$ and
$\beta\in\Omega$. If $z$ is a central element of $\l_\beta$, then
$\pr_\beta(z)$ is a central element of
$\varovee_\alpha\l_\alpha$.\end{theorem}

\begin{proof}
We denote $\pr_\beta(z)$ by $Z$. Hence $Z^c=\pr_\beta(z^c)$. From
Lemma \ref{LemmaCentralElements}, $z^c:=\s_\beta\backslash
z\in\l_\beta$. Therefore, by Axiom P2,
$Z^c\in\varovee_\alpha\l_\alpha$.

We now prove that $(Z,Z^c)M$. The proof for $(Z^c,Z)M$ is similar.
Let $R\in\varovee_\alpha\l_\alpha$ with $R\subseteq Z^c$. Let
$\mathbf{\s}=\prod_\alpha\s_\alpha$, $p\in \mathbf{\s}$,
$X=R\cup Z$, and $\alpha\in\Omega$. If $\alpha\ne \beta$, then
$X_\alpha[p]=\s_\alpha$ if $p_\beta\in z$ and
$X_\alpha[p]=R_\alpha[p]$ otherwise. On the other hand,
$X_\beta[p]=z\cup R_\beta[p]$, hence $X_\beta[p]\in\l_\beta$ by
Corollary \ref{CorollaryPatula}. As a consequence, $R\vee
Z=R\cup Z$, therefore $ (R\vee Z)\cap Z^c=(R\cup
Z)\cap Z^c=R$.
\end{proof}

\begin{theorem}[\cite{Ischi:2000}, Theorem 1]
Let $\{\l_\alpha\}_{\alpha\in\Omega}$ be a family of
orthocomplemented simple closure spaces on $\s_\alpha$. Suppose
that one of the following assumptions holds.
\begin{enumerate}
\item $\Omega$ is finite.
\item For all $\alpha\in\Omega$, $\l_\alpha$ has the covering
property, and for all $p\ne q\in\s_\alpha$ having the same central
cover, $p\vee q$ contains an infinite number of atoms.
\end{enumerate}
Then, $\varowedge_\alpha\l_\alpha$ is irreducible if and if all
$\l_\alpha$'s are irreducible.
\end{theorem}
\section{Automorphisms}\label{SectionAutomorphisms}

In this section we prove the following result. Let
$\l_1,\cdots,\l_n$ be simple closure spaces different from $2$,
$\l\in\sep(\l_1,\cdots,\l_n)$ and $u\in\cCal(\l,\l)$ sending atoms
to atoms. If $u$ is {\it large} (see Definition
\ref{DefinitionLarge} below), then there is a permutation $f$ of
$\{1,\cdots,n\}$ and arrows $v_i\in\cCal(\l_i,\l_{f(i)})$ sending
atoms to atoms such that for any atom $p$ of $\l$,
$u(p)_{f(i)}=v_i(p_i)$. We need some hypotheses on each $\l_i$
which are true for instance if each $\l_i$ is irreducible
orthocomplemented with the covering property or an irreducible
DAC-lattice. Note that our hypotheses imply irreducibility. Note
also that if $u$ is an automorphism, then $u$ is large.

\begin{definition}\label{DefinitionConnected}
Let $\l$ be a simple closure space on $\s$. We say that $\l$ is
{\it weakly connected} if $\l\ne 2$ and if there is a {\it
connected covering} of $\s$, that is a family of subsets
$\{A^\gamma\subseteq\s\,;\,\gamma\in\sigma\}$ such that
\begin{enumerate}
\item $\s=\bigcup \{A^\gamma\,;\,\gamma\in \sigma\}$ and $\abs{
A^\gamma}\geq 2$ for all $\gamma\in\sigma$,
\item for all $\gamma\in\sigma$ and all $p\ne q\in A^\gamma$,
$p\vee q$ contains a third atom,
\item for all $p,\,q\in\s$, there is a finite subset
$\{\gamma_1,\cdots,\gamma_n\}\subseteq\sigma$ such that $p\in
A^{\gamma_1}$, $q\in A^{\gamma_n}$, and such that
$\abs{A^{\gamma_i}\cap A^{\gamma_{i+1}}} \geq 2$ for all $1\leq
i\leq n-1$.
\end{enumerate}
We say that $\l$ is {\it connected} if $\l\ne 2$ and for all
$p,\,q\in\s$, $p\vee q$ contains a third atom, say $r$, such
that $p\in q\vee r$ and $q\in p\vee r$.
\end{definition}

\begin{remark}
Note that in part (2) of Definition \ref{DefinitionConnected}, it is
not required that the third atom under $p\vee q$ is in
$A^\gamma$. Note also that by Corollary \ref{CorollaryPatula},
weakly connected implies irreducible. Finally, let $\l$ be a
simple closure space. Then, if $\l\ne 2$ and $\l$ is irreducible
orthocomplemented with the covering property or an irreducible
DAC-lattice, then $\l$ is connected.
\end{remark}

\begin{definition}\label{DefinitionLarge}
Let $\{\l_\alpha\}_{\alpha\in\Omega}$ be a family of simple
closure spaces on $\s_\alpha$ different from $2$. Let
$\l\in\sep(\l_\alpha,\alpha\in\Omega)$,
$\mathbf{\s}=\prod_\alpha\s_\alpha$, and let $u\in\cCal(\l,\l)$.
We say that $u$ is {\it large} if for all $\beta\in\Omega$ and
$p\in\mathbf{\s}$, $u(p[\s_\beta])$ is not an atom of $\l$, and
$u(1)\nsubseteq \pr_\beta(p_\beta)$.\end{definition}

\begin{lemma}\label{LemmaJoinP1diffQ1andP2diffQ2}
Let $\{\l_\alpha\}_{\alpha\in\Omega}$ be a family of simple
closure spaces on $\s_\alpha$,
$\mathbf{\s}=\prod_\alpha\s_\alpha$, and let $\l$ be a simple
closure space on $\mathbf{\s}$. Suppose that Axiom P2 holds in
$\l$. Let $p,\,q\in\mathbf{\s}$.
\begin{enumerate}
\item If $p_\beta\ne q_\beta$ for at least two $\beta\in\Omega$,
then $p\vee q=p\cup q$.
\item For all $\beta\ne \gamma\in\Omega$ and for all
$b\in\l_\beta$ and $c\in\l_\gamma$ such that $p_\beta\in b$ and
$p_\gamma\in c$, $p[b,\beta]\vee p[c,\gamma]=p[b,\beta]\cup
p[c,\gamma]$.
\end{enumerate}
\end{lemma}

\begin{proof} (1) Let $r\in p\vee q$, and suppose that
$p_\beta\ne q_\beta$ and $p_\gamma\ne q_\gamma$ for some $\beta\ne
\gamma\in\Omega$. We must show that $r\in p\cup q$, {\it i.e.},
that $r=p$ or $r=q$. Now
\[\begin{split}
p,\,q\in\lp\pr_\beta(p_\beta)\cup
\pr_\gamma(q_\gamma)\rp&\cap\lp\pr_\beta(q_\beta)\cup
\pr_\gamma(p_\gamma)\rp\\
&=\lp\pr_\beta(p_\beta)\cap \pr_\gamma(p_\gamma)\rp\cup
\lp\pr_\beta(q_\beta)\cap \pr_\gamma(q_\gamma)\rp\,
,\end{split}\]
because
\[\pr_\beta(p_\beta)\cap\pr_\beta(q_\beta)=\emptyset=
\pr_\gamma(p_\gamma)\cap\pr_\gamma(q_\gamma) \, ,\]
since by hypothesis, $p_\beta\ne q_\beta$ and $p_\gamma\ne
q_\gamma$.

Now, by Axiom P2, $\pr_\alpha(x_\alpha)\in \l$ for all
$x\in\prod_\alpha\l_\alpha$ and all $\alpha\in\Omega$. Therefore,
\[ \lp\pr_\beta(p_\beta)\cup
\pr_\gamma(q_\gamma)\rp\cap\lp\pr_\beta(q_\beta)\cup
\pr_\gamma(p_\gamma)\rp\in \l .\]
As a consequence,
\[
p\vee q\subseteq \lp\pr_\beta(p_\beta)\cap
\pr_\gamma(p_\gamma)\rp\cup \lp\pr_\beta(q_\beta)\cap
\pr_\gamma(q_\gamma)\rp .\]
It follows that either $r_\beta=p_\beta$ and $r_\gamma=p_\gamma$,
or $r_\beta=q_\beta$ and $r_\gamma=q_\gamma$. Since this holds for
every pair of indices $\beta$, $\gamma$ at which $p$ and $q$
differ, either $r=p$ or $r=q$.

(2) Let $\beta\ne\gamma\in\Omega$, $b\in\l_\beta$,
$c\in\l_\gamma$, and $p,\, q\in\mathbf{\s}$ such that $p_\beta\in
b$ and $p_\gamma\in c$. We must show that $p[c,\gamma]\cup
p[b,\beta]\in\l$. Now, since $p_\beta\in b$ and $p_\gamma\in c$,
we have
\[\bigcap_{\alpha\ne\beta,\gamma}\pr_\beta(b)\cap
\pr_\gamma(c)\cap(\pr_\beta(p_\beta)\cup
\pr_\gamma(p_\gamma))
\cap\pr_\alpha(p_\alpha)=p[c,\gamma]\cup
p[b,\beta] .\]
To conclude, it suffices to note that by Axiom P2, the subsets of
$\mathbf{\s}$ $\pr_\beta(b)$, $\pr_\gamma(c)$,
$\pr_\beta(p_\beta)\cup \pr_\gamma(p_\gamma)$, and
$\pr_\alpha(p_\alpha)$ are elements of $\l$.
\end{proof}

\begin{theorem}\label{TheoremFactMaps}
Let $\Omega$  be a finite set and $\{\l_i\}_{i\in\Omega}$ a finite
family of connected simple closure spaces on $\s_i$. Let
$\l\in\sep(\l_i,i\in\Omega)$ and $u\in\cCal(\l,\l)$ large, sending
atoms to atoms. Then there is a bijection $f$ of $\Omega$, and for
each $i\in\Omega$, there is $v_i\in \cCal(\l_i,\l_{f(i)})$ sending
atoms to atoms such that $u(p)_{f(i)}=v_i(p_i)$ for all
$p\in\mathbf{\s}$ and $i\in\Omega$.
\end{theorem}

\begin{proof} The proof is similar to the proof of Theorem 3 in
\cite{Ischi:2000}.

(1) Let $p\in\mathbf{\s}$ and $j\in\Omega$. {\it Claim}: There is
$k \in\Omega$ such that
\[u(p[\s_j])\subseteq u(p)[\s_k] .\]
[{\it Proof}\nobreak : Since $\l_j$ is connected, for all
$q_j\in\s_j$ different from $p_j$, $p_j\vee q_j$ contains a
third atom, say $r_j$, and $p_j\in q_j\vee r_j$ and $q_j\in
p_j\vee r_j $. Suppose that $u(p)_k\ne u(p[q_j])_k$ for at least
two indices $k$. Then, by Lemma \ref{LemmaJoinP1diffQ1andP2diffQ2}
part (1) and Lemma \ref{LemmaAlex}.3,
\[u(p[r_j])\in u(p\vee p[q_j])=u(p)\vee u(p[q_j])=
u(p)\cup u(p[q_j]) .\]
Assume for instance that $u(p[r_j])=u(p)$. Then,
\[u(p[q_j])\in u(p\vee p[r_j])=u(p)\vee u(p[r_j])=u(p)\, ,\]
a contradiction, which proves the claim.]

(2) Let $j\in\Omega$. Denote the $k$ of part (1) by $f(j,p)$. {\it
Claim}: The map $p\mapsto f(j,p)$ is constant. [{\it
Proof}\nobreak : Let $p,\,q\in\mathbf{\s}$ differ only by one
component, say $j'\ne j$, {\it i.e.}, $p_{j'}\ne q_{j'}$, and
$p_i=q_i$, for all $i\ne j'$. Suppose that $f(j,p)\ne f(j,q)$.
Write $k:=f(j,p)$ and $k':=f(j,q)$.

(2.1) We first prove that
\[u(p[\s_j ])\vee u(q[\s_j])=u(p[\s_j])\cup u(q[\s_j]) .\]
By hypothesis,
\[u(p[\s_j])\subseteq u(p)[\s_k]\ \mbox{and}\ u(q[\s_j])\subseteq
u(q)[\s_{k'}] .\]
Hence, for all $r_j\in\s_j$, we have
\[u(p)_l=u(p[r_j])_l\,\forall l\ne k\ \mbox{and}\ u(q)_m=
u(q[r_j])_m\,\forall
m\ne k' .\]

Since $u$ is large, $u(p[\s_j])$ is not an atom, hence there is
$r_j\in\s_j$ such that
\[u(p[r_j])_k\ne u(q)_k .\]
Note that $u(q)_k=u(q[r_j])_k$. Therefore, since $p[r_j]$ and
$q[r_j]$ differ only by one component (namely $j'$), by part (1) we
have
\[u(p[r_j])_l=u(q[r_j])_l\,\forall l\ne k .\]
As a consequence,
\[u(q)_l=u(q[r_j])_l=u(p[r_j])_l=u(p)_l\,\forall l\ne k,\,k'\, ,\]
and
\[u(p)_{k'}=u(p[r_j])_{k'}=u(q[r_j])_{k'}\in\pi_{k'}(u(q[\s_j]))\,
.\]

On the other hand, since $u(q[\s_j])$ is not an atom, there is
$s_j\in\s_j$ such that
\[u(q[s_j])_{k'}\ne u(p)_{k'} .\]
Note that $u(p)_{k'}=u(p[s_j])_{k'}$. Therefore, since $q[s_j]$
and $p[s_j]$ differ only by on component (namely $j'$), by part (1)
we have
\[u(q[s_j])_l=u(p[s_j])_l\,\forall l\ne k' .\]
As a consequence,
\[u(q)_k=u(q[s_j])_k=u(p[s_j])_k\in\pi_k(u(p[\s_j])) .\]

To summarize, we have proved that $u(p)_l=u(q)_l$ for all $l\ne
k,\,k'$, $u(p)_{k'}\in\pi_{k'}(u(q[\s_j]))$, and $u(q)_k\in
\pi_k(u(p[\s_j]))$. As a consequence, the statement follows from
Lemma \ref{LemmaJoinP1diffQ1andP2diffQ2} part (2).

(2.2) Since $\l_{j'}$ is connected, there is $s_{j'}\in
p_{j'}\vee q_{j'}$ such that $p_{j'}\in q_{j'}\vee s_{j'}$
and $q_{j'}\in p_{j'}\vee s_{j'}$. Let $r=p[s_{j'}]$. Then, by
Lemma \ref{LemmaAlex}.3,
\[r[\s_j]\subseteq p[\s_j]\vee q[\s_j]\, ,\hspace{0.5cm}
p[\s_j]\subseteq q[\s_j]\vee r[\s_j]\hspace{0.3cm}
\mbox{and}\hspace{0.3cm} q[\s_j]\subseteq p[\s_j]\vee r[\s_j]\,
.\]
Now, by part (2.1),
\[u(p[\s_j]\vee q[\s_j])=u(p[\s_j])\vee u(q[\s_j])=u(p[\s_j])
\cup u(q[\s_j]) .\]
As a consequence, by part (1), $u(r[\s_j])\subseteq u(p[\s_j])$ or
$u(r[\s_j])\subseteq u(q[\s_j])$. Assume for instance that
$u(r[\s_j])\subseteq u(p[\s_j])$. Then,
\[u(q[\s_j])\subseteq  u(p[\s_j ]\vee r[\s_j])=u(p[\s_j])\vee
u(r[\s_j])=u(p[\s_j])\, ,\]
a contradiction. Hence, we have proved that if
$p,\,q\in\mathbf{\s}$ differ only by one component, then
$f(j,p)=f(j,q)$.

(2.3) Suppose now that $p$ and $q$ differ by more than one
component. Since $\Omega$ is finite, there is
$s_1,\cdots,s_n\in\mathbf{\s}$ such that $s^1=q$, $s^n=p$, and
such that for all $1\leq i\leq n-1$, $s^i$ and $s^{i+1}$ differ
only by one component. Therefore,
\[f(j,q)=f(j,s^1)=f(j,s^2)=\cdots=f(j,s^n)=f(j,p)\, ,\]
and we are done.]

(3) Let $p_0\in\mathbf{\s}$. Define $f:\Omega\rightarrow\Omega$ as
$f(i):=f(i,p_0)$. Note that by part (2), $f$ does not depend on the
choice of $p_0$. {\it Claim}: The map $f$ is surjective. [{\it
Proof}\nobreak : Let $k\in\Omega$. Suppose that for all
$i\in\Omega$ $f(i)\ne k$. Let $p,\, q\in\mathbf{\s}$ that differ
only by one component, say $j$. Then $p[\s_j]=q[\s_j]$, therefore
$u(p[\s_j])=u(q[\s_j])$. Moreover, by part (1), there is $k'\ne k$
such $u(p[\s_j ])\subseteq u(p)[\s_{k'}]$. As a consequence,
$u(p)_k=u(q)_k$.

Since $\Omega$ is finite, by the same argument as in part (2.3), we
find that $u(p)_k=u(q)_k$, for all $p,\,q\in\mathbf{\s}$. As a
consequence, $u(\mathbf{\s})\subseteq \pr_k(u(p_0)_k)$, a
contradiction since $u$ is large.]

(4) Let $p_0\in\mathbf{\s}$ and $j\in\Omega$. Define $v_j :\l_j
\rightarrow\l_{f(j )}$ as
\[v_j (a_j ):=\pi_{f(j )}(u(p_0[a_j,j ])) .\]
{\it Claim}: $v_j $ does not depend on the choice of $p_0$. [{\it
Proof}\nobreak : Let $q\in\mathbf{\s}$ that differs from $p_0$
only by one component, say $j'\ne j $. Then, by Lemma
\ref{LemmaAlex}.3, we have
\[\begin{split}\pi_{f(j )}(u(q[a_j ]))&=
\pi_{f(j)}(u(\vee \{q[r_j]\,;\,r_j\in
a_j\}))\\
&=\pi_{f(j)}(\vee \{u(q[r_j])\,;\,r_j\in a_j\})=\bigvee_{r_j\in
a_j} u(q[r_j])_{f(j)}\, ,\end{split}\]
and the same formula holds for $\pi_{f(j)}(u(p_0[a_j]))$. Now
\[u(q[r_j ])_{f(j )}=\pi_{f(j )}(u(s[\s_{j'}]))=u(p_0[r_j ])_{f(j
)}\, ,\]
where $s=q[r_j]$. As a consequence,
$\pi_{f(j)}(u(q[a_j]))=\pi_{f(j)}(u(p_0[a_j]))$.

Since $\Omega$ is finite, by the same argument as in part (2.3), we
find that 
\[\pi_{f(j )}(u(q[a_j ]))=\pi_{f(j )}(u(p_0[a_j ]))\]
for all $q\in\mathbf{\s}$.]

It remains to check that $v_j$ preserves arbitrary joins. Let
$\omega\subseteq\l_j$. Then, by Lemma \ref{LemmaAlex}.3,
\[\begin{split}
v_j(\vee\omega)&=\pi_{f(j)}(u(p_0[\vee\omega]))=
\pi_{f(j)}(u(\vee\{p_0[x,j]\,;\,x
\in\omega\}))\\
&=\pi_{f(j)}(\vee\{u(p_0[x,j])\,;\,x\in\omega\})\\
&=\bigvee_{x\in\omega}\pi_{f(j)}(u(p_0[x,j]))=
\bigvee_{x\in\omega}v_j(x) .\end{split}\]
\end{proof}

\begin{corollary}
If the $u$ in Theorem \ref{TheoremFactMaps} is an automorphism,
then all $v_i$'s are isomorphisms.\end{corollary}

\begin{theorem}\label{TheoremAutoFactorwithWeaklyConnected} If the
$u$ in Theorem \ref{TheoremFactMaps} is an automorphism, the
statement remains true if for all $i\in\Omega$, $\l_i$ is weakly
connected.\end{theorem}

\begin{proof}
The proof is similar as in Theorem \ref{TheoremFactMaps}. We only
sketch the arguments that must be modified. Note that since $u$ is
an automorphism and $\l_i\ne 2$ for all $i\in\Omega$, $u$ is
large.

(1) Since $\l_j$ is weakly connected, there is
$\gamma_0\in\sigma^j$ such that $p_j\in A_j^{\gamma_0}$. By
hypothesis and Lemma \ref{LemmaAlex}.3, for all $q_j\in
A_j^{\gamma_0}$, $p\vee (p[q_j])$ contains a third atom, hence
also $u(p)\vee u(p[q_j])$ since $u$ is injective. As a
consequence, there is $k_{\gamma_0}\in\Omega$ such that
$u(p[A_j^{\gamma_0} ])\subseteq u(p)[\s_{k_{\gamma_0}}]$.
Moreover, since $u$ is injective, by the third hypothesis in
Definition \ref{DefinitionConnected}, the map $\gamma\mapsto
k_\gamma$ is constant. Therefore, since $\cup A^\gamma_j=\s_j$,
we find that $u(p[\s_j ])\subseteq u(p)[\s_k ]$.

(2) Take $p,\,q\in\mathbf{\s}$ that differ only by one component
such that $q_{j'}$ and $p_{j'}$ are in the same $A_j^\gamma$.

(2.2) By hypothesis, $p_{j'}\vee q_{j'}$ contains a third atom,
say $r_{j'}$, therefore $r[\s_j]\subseteq p[\s_j]\vee q[\s_j]$,
hence, $u(r[\s_j])=u(p[\s_j])$ or $u(r[\s_j])=u(q[\s_j])$, a
contradiction since $u$ is injective. As a consequence,
$f(j,q)=f(j,p)$. Now, since $u$ is injective, by the third
hypothesis in definition \ref{DefinitionConnected}, we find that
$f(j,q)=f(j,p)$, for all $p,\,q\in\mathbf{\s}$ that differ only by
the component $j'$.
\end{proof}
\section*{Acknowledgments}
A part of this work was done during a stay at McGill. In this
connection, I would like to thanks M. Barr for his hospitality.
\bibliographystyle{abbrv}
\bibliography{../../Bibliographie/References}
\end{document}